\documentclass[submission]{dmtcs}

\usepackage[latin1]{inputenc}
\usepackage{subfigure}
\usepackage[round]{natbib}

\RCSdef$Revision: 1.3 $\endRCSdef
\rcsMajMin

\newtheorem{theorem}{Theorem}

\newtheorem{corollary}[theorem]{Corollary}

\newtheorem{definition}[theorem]{Definition}

\newtheorem{lemma}[theorem]{Lemma}
\newtheorem{notation}[theorem]{Notation}

\newtheorem{remark}[theorem]{Remark}

\newcommand {\G}{\Gamma}

\newcommand {\la}{\lambda}

\def\b{\beta}

\def\G{\Gamma}

\def\l{\lambda}

\def\s{\sigma}

\def\f{\rightarrow}
\def\tr{\triangleright}

\def\v{\vdash}

\def\ou{\vee}
\def\et{\wedge}
\def\<{\langle}
\def\>{\rangle}
\def\F{\displaystyle\frac}

\author{Karim Nour\addressmark{1}\thanks{Karim.Nour@univ-savoie.fr}}

\title[Classical Combinatory Logic]{Classical Combinatory Logic}

\address{\addressmark{1} LAMA - Equipe de logique , Universit\'e de
Savoie ,  F-73376 Le Bourget du Lac, France
}

\revision{\rcsMaj}

\received{...}
\revised{...}
\accepted{...}

\begin{document}

\maketitle
\begin{abstract} 
Combinatory logic shows that bound variables can be eliminated
without loss of expressiveness. It has applications both in the
foundations of mathematics and in the implementation of functional
programming languages. The original combinatory calculus
corresponds to minimal implicative logic written in a system ``\`a
la Hilbert''. We present in this paper a combinatory logic which
corresponds to propositional classical logic. This system is
equivalent to the system $\l^{Sym}_{Prop}$ of Barbanera and Berardi.
\end{abstract}

\clearpage
\tableofcontents

\section{Introduction}

Combinatory logic started with a paper by \cite{Sch}. The
aim was an elimination of bound variables.  He proved that it is
possible to reduce the logic to a language consisting of one
constructor (the application) and some primitive constants. This work
was continued by \cite{Curry1} who introduced the syntax of the terms of
combinatory logic. At about the same time, \cite{Chu} introduced the
lambda-calculus as a new way to study the concept of rule. Originally
his purpose was to provide a foundation for mathematics. Combinatory
logic and lambda-calculus, in their type-free version,
generate essentially the same algebraic and logic structures. The
original combinatory calculus corresponds to minimal implicative logic
presented in a system ``\`a la Hilbert''. The codings between
combinatory logic and simply typed calculus preserve types. Research
on combinatory logic has been continued essentially by Curry's
students, \cite{Hin2}.

Since it has been understood that the Curry-Howard isomorphism
relating proofs and programs can be extended to classical logic,
various systems have been introduced: the $\l_c$-calculus
(\cite{Kri}), the $\la_{exn}$-calculus (\cite{deG4}), the $\l
\mu$-calculus (\cite{Par1}), the $\lambda^{Sym}$-calculus
(\cite{BaBe}), the $\lambda_{\Delta}$-calculus (\cite{ReSo}), 
the $\overline{\lambda}\mu\tilde{\mu}$-calculus
(\cite{cuhe}), the dual calculus (\cite{wad1}) ... 
All these calculi are based on logical systems presented either  
in natural deduction or in sequent calculus. 

We wish to define a combinatory calculus which corresponds to classical
logic presented ``\`a la Hilbert''. There are two ways to define such a
calculus:

- Add new combinators for the axioms which define classical logic over
minimal logic and give the corresponding reduction rules.

- Code by combinators an existing calculus based on classical logic.

The first way gives a very ``artificial'' solution. The reduction
rules for the new combinators are rather complicated. For the second
way, it is necessary to choose a system such that the reduction
rules erase the abstractors (i.e. the right-hant side of the reduction
rules should not introduce new abstractions). One of these calculi is the
$\lambda^{Sym}$-calculus of Barbanera and Berardi.

We present in this paper the  $\l^{Sym}$-calculus and the new
combinatory calculus CCL. We also explain how to encode each calculus into
the other.

The paper is organized as follows. In section 2, we give the
syntax of the terms and the reduction rules of the system
$\l^{Sym}_{Prop}$.  We introduce, in section 3, the
syntax of the terms and the reduction rules of the system
CCL. We encode, in section 4, the system $\l^{Sym}_{Prop}$ into the system CCL
and we encode, in section 5, the system CCL into the system
$\l^{Sym}_{Prop}$. We conclude with some future work.

\section{The system $\l^{Sym}_{Prop}$}

\begin{definition}
\begin{enumerate}
\item We have two sets of base types ${\cal A} = \{a,b,...\}$
and  ${\cal A}^{\perp} = \{a^{\perp},b^{\perp},...\}$.

\item The set of $m$-types is defined by the following grammar:

$$A ::= {\cal A} \mid {\cal A}^{\perp} \mid A \et A \mid A \ou A$$

\item The set of types is defined by the following grammar:

$$C ::= A \mid \perp $$

\item We define the negation $A^{\perp}$ of an $m$-type as follows:
\begin{itemize}
\item $(a)^{\perp} = a^{\perp}$
\item $(a^{\perp})^{\perp} = a$
\item $(A \et B)^{\perp} = A^{\perp} \ou B^{\perp}$
\item $(A \ou B)^{\perp} = A^{\perp} \et B^{\perp}$
\end{itemize}
\end{enumerate}
\end{definition}

\begin{lemma}
For all $m$-type $A$, $A ^{\perp \perp} = A$.
\end{lemma}

\begin{proof}
By induction on $A$.
\end{proof}

\begin{definition}
\begin{enumerate}
\item The terms of the system $\l^{Sym}_{Prop}$ (called $\l_s$-terms) are defined (in the
natural deduction style) by the following rules:
\begin{center}
$\F{}{\G , x : A \v x : A}$\\[0.4cm]
$\F{\G \v u : A \;\;\; \G \v v : B}{\G \v \<u,v\> : A \et B}$
$\;\;\;\;\;\;\;\;\;\;\;\;$ $\F{\G \v t : A}{\G \v \s_1(t) : A \ou
B}$ $\;\;\;\;\;\;\;\;\;\;\;\;$ $\F{\G \v t : B}{\G \v \s_2(t) : A
\ou B}$
\\[0.4cm]
$\F{\G , x : A \v t : \perp}{\G \v \l x.t : A^{\perp}}$
$\;\;\;\;\;\;\;\;\;\;\;\;$ $\F{\G \v u : A^{\perp} \;\;\; \G \v v
: A}{\G \v u \star v : \perp}$
\end{center}
We write $\G \v_{\l_s} t : A$, if we can type the $\l_s$-term
$t$ by the type $A$ using the set of declaration of variables $\G$.
\item The reduction rules are the following:
\begin{center}
\begin{tabular}{ccc}
$(\l x.u) \star v$ & $\f_{\b}$ & $u [ x:= v]$\\
$v \star (\l x.u)$ & $\f_{\b^{\perp}}$ & $u [ x:= v]$\\
$\l x. (u \star x)$ & $\f_{\eta}$ & $\;\;\;\;\;\;\;\;\;$ $u$ $\;\;\;$ {\rm (1)}\\
$\l x. (x \star u)$ & $\f_{\eta^{\perp}}$ & $\;\;\;\;\;\;\;\;\;$ $u$ $\;\;\;$ {\rm  (1)}\\
$\<u , v\> \star \s_1(w)$ & $\f_{\pi_1}$ & $u \star w$\\
$\<u , v\> \star \s_2(w)$ & $\f_{\pi_2}$ & $v \star w$\\
$\s_1(w) \star \<u , v\>$ & $\f_{\pi^{\perp}_1}$ & $w \star u$\\
$\s_2(w) \star \<u , v\>$ & $\f_{\pi^{\perp}_2}$ & $w \star v$\\
$u[x:=v]$ & $\f_{triv}$ & $\;\;\;\;\;\;\;\;\;$ $v$ $\;\;\;$ {\rm (2)}
\end{tabular}
\end{center}
{\rm (1)} if $x \not \in Fv(u)$

{\rm (2)} if $u$ and $v$ are $\l_s$-terms with type $\perp$, $x$
occurs only one time in $u$ and $u \neq x$. In this case $v = v_1
\star v_2$ and $\l y.x$ is a sub-term of $u$.

\item  We denote by $\f$ the one of previous rules. The
 transitive (resp. reflexive and transitive) closure of $\f$ is
 denoted by $\f^+$ (resp. $\f^*$).

\item We denote the $\l_s$-terms by small letters like $t,u,v,...$.
\end{enumerate}
\end{definition}

\begin{remark}
The reduction $\f^*$ is not confluent. For example $(\l x.(y \star z))
\star (\l x'.(y' \star z'))$ reduces both to $y \star z$ and to $y' \star z'$.
\end{remark}

\begin{theorem}[Subject reduction]
If $\G \v_{\l_s} u : A$ and $u \f^* v$, then $\G \v_{\l_s} v : A$.
\end{theorem}

\begin{proof}
It is enough to check that every reduction rule
preseves the type.
\end{proof}

\begin{theorem}[Strong normalization]\label{snls}
Every $\l_s$-term is strongly normalizing.
\end{theorem}

\begin{proof}
See \cite{BaBe}.
\end{proof}

\begin{remark}
\cite{BaBe} proved the strong normalization
of the $\l^{Sym}_{Prop}$-calculus by using candidates of reducibility
but, unlike the usual construction (for example for Girard's system
$F$), the definition of the interpretation of a type needs a rather
complex fix-point operation. This proof is highly non
arithmetical. P. Battyanyi recently gave an arithmetical proof of this
result by using the methods developed in \cite{tlca} to show the
strong normalization of systems $\l \mu \mu'$- calculus and
$\overline{\l} \mu \tilde{\mu}$-calculus.
\end{remark}

\section{The system CCL}

\begin{definition}
\begin{enumerate}
\item We use the same types as in section 2. The
terms of the system {\rm CCL} (called $c$-terms) are defined (in the
Hilbert style) by the following rules:
\begin{center}
$\F{}{\G , x : A \v x : A}$\\[0.4cm]

$\F{}{\G \v {\bf K} : A^{\perp} \ou (B \ou A)}$\\[0.4cm]

$\F{}{\G \v {\bf S} : (A \et (B \et C^{\perp})) \ou ((A \et B^{\perp}) \ou (A^{\perp} \ou C))}$\\[0.4cm]

$\F{}{\G \v {\bf C} : (A \et B) \ou ((A \et B^{\perp}) \ou A^{\perp})}$\\[0.4cm]

$\F{}{\G \v {\bf P} : A^{\perp} \ou (B^{\perp} \ou (A \et B))}$\\[0.4cm]

$\F{}{\G \v {\bf Q_1} : A^{\perp} \ou (A \ou B)}$
$\;\;\;\;\;\;\;\;\;\;\;\;$
$\F{}{\G \v {\bf Q_2} : B^{\perp} \ou (A \ou B)}$\\[0.4cm]

$\F{\G \v U : A^{\perp} \ou B \;\;\; \G \v V : A}{\G \v (U \; V) :
B}$ $\;\;\;\;\;\;\;\;\;\;\;\;$ $\F{\G \v U : A^{\perp} \;\;\;\G \v
V : A}{\G \v U \star V : \perp}$
\end{center}
Note that the typed rules does not change the set of declaration of variables.
We write $\G \v_{c} T : A$, if we can type the $c$-term $U$
by the type $A$ using the  set a declaration of variables $\G$.
\item Let $U, U_1,U_2, ...,U_n$ be $c$-terms. We write $(U \; U_1 \; U_2 \; ...
\; U_n)$\\
instead of $(...((U \; U_1) \; U_2) \; ... \; U_n)$.
\item The reduction rules are the following:
\begin{center}
\begin{tabular}{ccc}
$({\bf K} \; U \; V)$ & $\tr_{K}$ & $U$\\
$({\bf S} \; U \; V \; W)$ & $\tr_{S}$ & $((U \; W) \; (V \; W))$\\
$({\bf C} \; U \; V) \star W$ & $\tr_{C_r}$ & $(U \; W) \star (V \; W)$\\
$W \star ({\bf C} \; U \; V)$ & $\tr_{C_l}$ & $(U \; W) \star (V \; W)$\\
$({\bf C} \; ({\bf K} \; U) \; {\bf I})$ & $\tr_{e_r}$ & 
$\;\;\;\;\;\;\;\;\;\;\;\;$ $U$  $\;\;\;\;$ {\rm
(3)}\\
$({\bf C} \; {\bf I} \; ({\bf K} \; U))$ & $\tr_{e_l}$ &
$\;\;\;\;\;\;\;\;\;\;\;\;$ $U$  $\;\;\;$ {\rm
(3)}\\
$({\bf P} \; U \; V) \star ({\bf  Q_1} \; W)$ & $\tr_{pq_1}$ & $U \star W$\\
$({\bf P} \; U \; V) \star ({\bf Q_2} \; W)$ & $\tr_{pq_2}$ & $V \star W$\\
$({\bf Q_1} \; W) \star ({\bf P} \; U \; V)$  & $\tr_{qp_1}$ & $W \star U$\\
$({\bf Q_2} \; W) \star ({\bf P} \; U \; V)$  & $\tr_{qp_2}$
& $W \star V$\\
$W[x:=({\bf C} \; ({\bf K} \; U) \; ({\bf K} \; V))]$ &
$\;\;\tr_{simp}$ & $\;\;\;\;\;\;\;\;\;$ $U \star V$ $\;\;\;$ {\rm (4)}
\end{tabular}
\end{center}
{\rm (3)} where ${\bf I} = ({\bf S} \; {\bf K} \; {\bf K})$.

{\rm (4)} if $W$ is a $c$-term with type $\perp$.
\item  We denote by $\tr$ the one of previous rules. The
 transitive (resp. reflexive and transitive) closure of $\tr$ is
 denoted by $\tr^+$ (resp. $\tr^*$).
\item We denote the $c$-terms by capital letters like
$T,U,V,...$.
\end{enumerate}
\end{definition}

\begin{remark}
\begin{enumerate}
\item  We have $\v_{C} {\bf I} : A^{\perp} \ou A$ and,
for all $c$-term $T$, $({\bf I} \; T) \tr^* T$.
\item The reduction $\tr^*$ is not confluent. For example $({\bf C} \; ({\bf
K} \; y) \; ({\bf K} \; z)) \star ({\bf C} \; ({\bf K} \; y') \; ({\bf
K} \; z'))$ reduces both to $y \star z$ and to $y' \star z'$.
\end{enumerate}
\end{remark}

\begin{theorem}[subject reduction]
If $\G \v_{c} U : A$ and $U \tr^* V$, then $\G \v_{c} V : A$.
\end{theorem}

\begin{proof}
It is enough to check that every reduction rule
preserves the type.
\end{proof}

\begin{definition}
\begin{enumerate}
\item A $c$-term is said to be pre-term iff it does not contain the symbol
$\star$.
\item  A $c$-term $T$ is said to be star-term iff $T = U \star V$
for some pre-terms $U,V$.
\end{enumerate}
\end{definition}

\begin{lemma}\label{pre-star}
\begin{enumerate}
\item If $A$ is an m-type and $\G \v_{c} T : A$, then $T$ is a pre-term.
\item If $\G \v_{c} T : \perp$, then $T$ is a star-term.
\end{enumerate}
\end{lemma}

\begin{proof}
Easy.
\end{proof}

\begin{corollary}
A $c$-term is either a pre-term or a star-term.
\end{corollary}

\begin{proof}
By lemma \ref{pre-star}.
\end{proof}

\section{The encoding of $\l^{Sym}_{Prop}$ into CCL}

\begin{definition}
The function $\phi : \l^{Sym}_{Prop} \f {\rm CCL}$ is defined as follows:
\begin{itemize}
\item $\phi(x) = x$
\item $\phi(\l x.t) = l_x(\phi(t))$
\item $\phi(u \star v) = \phi(u) \star \phi(v)$
\item $\phi(\<u,v\>) = ({\bf P} \; \phi(u) \; \phi(v))$
\item $\phi(\s_1(t)) = ({\bf Q_1} \;  \phi(t))$
\item $\phi(\s_2(t)) = ({\bf Q_2} \; \phi(t))$
\end{itemize}
where
\begin{itemize}
\item $l_x(x) = {\bf I}$
\item $l_x(T) = ({\bf K} \; T)$ if $T$ is a pre-term and $x \not \in Var(T)$
\item $l_x((U \; V)) = ({\bf S} \;\; l_x(U) \; l_x(V))$ if $x \in
Var((U \; V))$
\item $l_x(U \star V) = ({\bf C} \;\; l_x(U) \; l_x(V))$
\end{itemize}
\end{definition}

\begin{lemma}\label{conctype}
Let $A$ and $B$ be m-types.
\begin{enumerate}
\item If $\G, x:A \v_{c} T: B$, then  $\G \v_{c} l_x(T):A^{\perp} \ou B$.
\item If $\G, x:A \v_{c} T:\perp$, then  $\G \v_{c} l_x(T):A^{\perp}$.
\end{enumerate}
\end{lemma}

\begin{proof}
1. By induction on $T$. 

2. Use $1.$
\end{proof}

\begin{theorem}
If $\G \v_{\l_s} t:A$, then  $\G \v_{c} \phi(t):A$.
\end{theorem}

\begin{proof}
By induction on the typing. Use lemma \ref{conctype}.
\end{proof}

\begin{lemma}\label{l_x}
\begin{enumerate}
\item If $U$ is a pre-term, then $(l_x(U) \, V) \tr^* U[x:=V]$.
\item If $U$ is a star-term, then $l_x(U) \star V \tr^* U[x:=V]$ and $V \star l_x(U) \tr^* U[x:=V]$.
\end{enumerate}
\end{lemma}

\begin{proof}
1. By induction on $U$. 

2. Use $1.$
\end{proof}

\begin{lemma}\label{phi}
\begin{enumerate}
\item If $V$ is a pre-term and $x \not \in Var(V)$, then $l_x(U[y:=V]) = l_x(U)[y:=V]$.
\item $\phi(u[y:=v]) = \phi(u)[y:=\phi(v)]$.
\end{enumerate}
\end{lemma}

\begin{proof}
1. By induction on $U$. 

2. By induction on $u$. Use $1.$
\end{proof}

\begin{remark}
As in $\l$-calculus, we do not have, in general, if $u \f v$, then
$\phi(u) \tr^+ \phi(v)$. The problem comes from the
${\beta}$-reductions ``under a lambda''.
\end{remark}

\begin{definition}
We write $u \f_{\omega} v$ if $v$ is obtained
by reducing in $u$ a redex which is not within the scope of a
$\l$-abstraction.
\end{definition}

\begin{theorem}
If $u \f_{\omega}  v$, then $\phi(u) \tr^+ \phi(v)$.
\end{theorem}

\begin{proof}
By induction on $u$. Use lemmas \ref{l_x} and \ref{phi}.
\end{proof}

\section{The encoding of CCL into $\l^{Sym}_{Prop}$}

\begin{notation}
Let $\pi_i t$ denote the $\l_s$-term $\l x. (t \star
\s_i(x))$ where $i \in \{1,2\}$ and $x \not
\in Fv(t)$. For each $i_1,...,i_n \in \{1,2\}$, let
$\pi_{i_1...i_n}t$ denote the $\l_s$-term $\pi_{i_1}...\pi_{i_n} t$.
\end{notation}

\begin{lemma}\label{pi}
\begin{enumerate}
\item $\pi_1 \<u,v\> \f^* u$ and $\pi_2 \<u,v\> \f^* v$.
\item If $\G \v_{\l_s} t : A \et B$, then $\G \v_{\l_s} \pi_1 t : A$
and $\G \v_{\l_s} \pi_2 t : B$.
\end{enumerate}
\end{lemma}

\begin{proof}
Easy.
\end{proof}

\begin{notation}
Let $[u,v]$ denote the $\l_s$-term $\l x. (u \star
\<v,x\>)$ where $x \not \in Fv(u) \cup Fv(v)$.
\end{notation}

\begin{lemma}\label{app}
\begin{enumerate}
\item $[\l x. u , v] \f^* \l y. u[x := \<v , z\>]$.
\item If $\G \v_{\l_s} u : A^{\perp} \ou B$ and $\G' \v_{\l_s} v : A$,
then $\G,\G' \v_{\l_s} [u,v] : B$.
\end{enumerate}
\end{lemma}

\begin{proof}
Easy.
\end{proof}

\begin{definition}
The function $\psi :  {\rm CCL} \f \l^{Sym}_{Prop}$ is defined as follows:
\begin{itemize}
\item $\psi(x) = x$
\item $\psi({\bf K}) = \l x. (\pi_1 x \star \pi_{22}x)$
\item $\psi({\bf S}) = \l x. ([[\pi_1 x , \pi_{122}x] , [\pi_{12}x ,\pi_{122}x ]] \star \pi_{222}x)$
\item $\psi({\bf C}) =  \l x. ([\pi_1 x , \pi_{22}x] \star [\pi_{12}x ,\pi_{22}x])$
\item $\psi({\bf P}) = \l x. (\<\pi_1 x , \pi_{12}x \> \star \pi_{22}x)$
\item $\psi({\bf Q_1}) = \l x. (\s_1(\pi_{1} x) \star \pi_{2} x)$
\item $\psi({\bf Q_2}) = \l x. (\s_2(\pi_{1} x) \star \pi_{2} x)$
\item $\psi((U \; V)) = [\psi(U),\psi(V)]$
\item $\psi(U \star V) = \psi(U) \star \psi(V)$
\end{itemize}
\end{definition}

\begin{theorem}
If $\G \v_{c} U : A$, then $\G \v_{\l_s} \psi(U) : A$.
\end{theorem}

\begin{proof}
Use lemmas \ref{pi} and \ref{app}.
\end{proof}

\begin{lemma}\label{subpsi}
$\psi(U[x:=V]) = \psi(U)[x:=\psi(V)]$.
\end{lemma}

\begin{proof}
By induction on $U$.
\end{proof}

\begin{theorem}\label{simcclls}
If $U \tr V$, then $\psi(U) \f^+ \psi(V)$.
\end{theorem}

\begin{proof}
The following are easy to check:

\begin{tabular}{ccc}
$[[\psi({\bf K}) , u] , v]$ & $\f^{+}$ & $u$\\
$[[[\psi({\bf S}) , u] , v ] , w]$ & $\f^{+}$ & $[[u , w ], [v , w]]$\\
$[\psi({\bf I}),u]$ & $\f^{+}$ & $u$\\
$[[\psi({\bf C}) , u] , v] \star w$ & $\f^{+}$ & $[u , w ] \star [v , w ]$\\
$w \star [[{\bf C} , u] , v]$ & $\f^{+}$ & $[u , w ] \star [v , w ]$\\
$[[\psi({\bf C}) , [\psi({\bf K}) , u ]] , \psi({\bf I})]$ & $\f^{+}$ & $u$\\
$[[\psi({\bf C}) ,\psi({\bf I})] , [\psi({\bf K}) , u]]$ & $\f^{+}$ & $u$\\
$[[\psi({\bf P}) , u] , v] \star [\psi({\bf Q_1}) , w]$ & $\f^{+}$ & $u \star w$\\
$[[\psi({\bf P}) , u] , v] \star [\psi({\bf Q_2}) , w]$ & $\f^{+}$ & $v \star w$\\
$[\psi({\bf Q_1}) , w] \star [[\psi({\bf P}) , u] , v]$  & $\f^{+}$ & $w \star u$\\
$[\psi({\bf Q_2}) , w] \star [[\psi({\bf P}) , u] , v]$ &
$\f^{+}$ & $w \star v$\\
$[[\psi({\bf C}) , [\psi({\bf K}) , u ]] , [\psi({\bf K}) , u ]]$
& $\f^{+}$ & $\l z. (u \star v)$
\end{tabular}

For the reduction rule $\tr_{simp}$, we use lemma \ref{subpsi}.
\end{proof}

\begin{theorem}[Strong normalization]\label{snccl}
Every $c$-term is strongly normalizing.
\end{theorem}

\begin{proof}
By theorems \ref{simcclls} and \ref{snls}.
\end{proof}

\section{Future work}

Although the strong normalization of the system CCL
follows from the one of the system $\l^{Sym}_{Prop}$ (see theorem
\ref{snccl}), R. David and I aim to prove directly 
this property. We wish to deduce a simpler proof of the strong
normalization of the system $\l^{Sym}_{Prop}$.  For that, it is
necessary to show a notion stronger than the strong normalization
because the coding, presented in section 4, does not simulate all
reductions. The verifications we made for the ordinary combinatory
logic are very promizing. 

In the original combinatory logic the reduction rules of {\bf K} and
{\bf S} do not allow ${\beta}$-reduction to be fully simulated (the
problem comes from the ${\beta}$-reductions ``under a
lambda''). Nevertheless, by adding an extensionality rule to
combinatory logic (i.e. $\forall \, x \, \{(F \, x) = (G \, x)\}
\Rightarrow F= G$) one obtains an equational theory that corresponds
exactly to $\beta\eta$-equivalence. The question is ``Is there
anything similar for CCL?''. This question is not an easy one because
CCL is not confluent. Consequently, a weaker notion than
extensionality would be needed.

\acknowledgements
\label{sec:ack}
I wish to thank Ren\'e David for helpful discussions.

\nocite{*}
\bibliographystyle{abbrvnat}
\bibliography{CCL}
\label{sec:biblio}

\end{document}